\documentclass{notices}
\usepackage{amsfonts,amssymb,amsmath,amscd}
\usepackage{graphicx,comment,colonequals}

\newcommand{\Z}{\mathbb{Z}}
\newcommand{\Q}{\mathbb{Q}}
\newcommand{\R}{\mathbb{R}}
\newcommand{\C}{\mathbb{C}}
\newcommand{\F}{\mathbb{F}}
\newcommand{\fm}{\mathfrak{m}}

\newcommand{\fq}{\mathfrak{q}}
\newcommand{\fa}{\mathfrak{a}}
\newcommand{\st}{\mathrm{st}}
\newtheorem{thm}{Theorem}
\newtheorem{defn}{Definition}

\newtheorem{comp}{Comparison}
\newtheorem{ex}{Example}

\title{
Differentiating by prime numbers
}

\author{
 Jack Jeffries
  \affil{
    The author is an assistant professor of mathematics at the University of Nebraska. His email address is jack.jeffries@unl.edu. His work is supported by NSF CAREER Award DMS-2044833.
    }
}

\begin{document}

\maketitle

\section*{}
It is likely a fair assumption that you, the reader, are not only familiar with but even quite adept at differentiating by $x$. What about differentiating by 13? That certainly didn't come up in my calculus class! From a calculus perspective, this is ridiculous: are we supposed to take a limit as 13 changes?

One notion of differentiating by 13, or any other prime number, is the notion of \emph{$p$-derivation} discovered independently by Joyal \cite{MR789309} and Buium \cite{MR1387233}. $p$-derivations have been put to use in a range of applications in algebra, number theory, and arithmetic geometry. Despite the wide range of sophisticated applications, and the fundamentally counterintuitive nature of the idea of differentiating by a number, $p$-derivations are elementary to define and inviting for exploration.

In this article, we will introduce $p$-derivations and give a few basic ways in which they really do act like derivatives by numbers; our hope is that you will be inspired and consider adding $p$-derivations to your own toolkit!

\section*{ $p$-derivations on $\Z$}
First we want to discuss differentiating one number $n$, by another, $p$; i.e., what we will call $p$-derivations on $\Z$. Before we succeed, we need to abandon the notion of derivative as a limit as the input varies by a small amount: the thing $p$ that we are differentiating by does not vary, and the thing $n$ that we are differentiating does not even have an input! Instead, we take a little inspiration from elementary number theory.

Let $p$ be a prime number. By Fermat's little theorem, for any integer $n$, we have
\[ n \equiv n^p \mod \ p,\]
so we can divide the difference $n-n^p$ by $p$. The starting point of our journey is that  not only \emph{can} we divide by $p$ here, but we \emph{should}.
The $p$-derivation on $\Z$ is the result of this process. Namely:

\begin{defn} For a prime number $p$, the \emph{$p$-derivation} on $\Z$ is defined as the function ${\delta_p:\Z\to \Z}$ given by the formula
\[ \delta_p(n) = \frac{n-n^p}{p}.\]
\end{defn}

So, in particular, there is the $2$-derivation on $\Z$ and the $13$-derivation  on $\Z$  given respectively by
\[ \delta_2(n) = \frac{n-n^2}{2} \ \ \text{and} \ \ \delta_{13}(n) = \frac{n-n^{13}}{13}.\]

Let's plug in a few values:
\begin{center}
\[ \begin{array}{ c | c c c }
n & \delta_2(n) & \delta_3(n) & \delta_5(n) \\
\hline
\vdots &\vdots &\vdots & \vdots \\ 
-4 & -10 & 20 & 204 \\
-3 & -6 & 8 & 48  \\
-2 & -3 & 2 & 6 \\
-1 & -1 & 0 & 0  \\
0 & 0 & 0 & 0 \\
1 & 0 & 0 & 0  \\
2 & -1 & -2 & -6 \\
3 & -3 & -8 &-48  \\
4 & -6 & -20 &-204 \\
5 & -10 & -40 &-624 \\
6 & -15 & -70 &-1554 \\
\vdots &\vdots &\vdots &\vdots 
\end{array}\]
\end{center}

A quick look at this table suggests a few observations, easily verified from the definition:
\begin{itemize}
\item Numbers are no longer ``constants'' in the sense of having derivative zero, but at least $0$ and $1$ are.
\item These functions are neither additive  nor multiplicative, e.g.:
\[ \delta_p(1) + \delta_p(1) \neq \delta_p(2),\]
 \[ \delta_p(1) \delta_p(2) \neq \delta_p(2).\]
\item $\delta_p$ is an odd function, at least for $p\neq 2$.
\item The outputs of $\delta_2$ are just the negatives of the triangular numbers.
\end{itemize}

We might also note that the outputs are very large in absolute value, and think that this operation is simply making a mess of our numbers. However, something more informative occurs if we think about largeness of the outputs from the point of view of $p$, namely, the \emph{$p$-adic order} of $n$---the number of copies of $p$ in its prime factorization. Writing $n=p^a m$ with $\gcd(m,p)=1$, if $a>0$, we get
\[ \delta_p(p^am) = \frac{p^a m - (p^a m)^p}{p} = p^{a-1} m (1-p^{ap-a}m^{p-1}).\]
Since $p\geq 2$ and $a\geq 1$, we must have $ap-a\geq 1$, so $p$ does not divide $1-p^{ap-a}m^{p-1}$. In particular, the $p$-derivation decreases the $p$-adic order of a multiple of $p$ by exactly one. This leads to our first comparison with old-fashioned $\frac{d}{dx}$:

\begin{comp}[Order-decreasing property]{\ } 
\begin{itemize}
\item If $f\in \R[x]$ is a polynomial and $x=r$ is a root of $f$ of multiplicity $a>0$, then $x=r$ is a root of the polynomial $\frac{d}{dx}(f(x))$ of multiplicity $a-1$.
\item If $n$ is an integer and $p$ is a prime factor of $n$ of multiplicity $a>0$, then $p$ is a prime factor of the integer $\delta_p(n)$ of multiplicity $a-1$.
\end{itemize}
\end{comp}
In particular, if $r$ is a simple root or $p$ is a simple factor, then it is no longer a root or factor of $\frac{d}{dx}(f(x))$ or $\delta_p(n)$ respectively. 

Let's check this against our table: the numbers $-2,2$, and $6$ that were divisible by $2$ but not $4$ result in odd numbers when we apply $\delta_2$, whereas $\pm 4$ returned even numbers no longer divisible by $4$. Note that this order-decreasing property says nothing about what happens when you apply $\delta_2$ to an odd number, and indeed, based on the table we observe that even and odd numbers can result. You can convince yourself that
\[ \delta_2(n) \ \text{is} \ \begin{cases} \text{even} \ \text{if} \ n\equiv 0,1 \mod 4 \\ \text{odd} \ \text{if} \ n\equiv 2,3 \mod 4.\end{cases}\]

We've observed already that these $p$-derivations on $\Z$ are not additive. This can be a bit unsettling for those of us (like myself) who are usually accustomed to the luxury of additive operators. However, any function satisfying the order-decreasing property of $\delta_p$ above must not be additive, since an additive function has to take multiples of $p$ to multiples of $p$. However, the error term can be made concrete:
\[\begin{aligned} \delta_p(m+n) -( \delta_p(m) + \delta_p(n) ) &= \frac{m^p + n^p -(m+n)^p}{p} \\&= -\sum_{i=1}^{p-1} \frac{\binom{p}{i}}{p} m^i n^{p-i}.\end{aligned}\]
All of the binomial coefficients $\binom{p}{i}$ appearing above are multiples of $p$, so this expression is, given a particular value of $p$, a particular polynomial in $m$ and $n$ with integer coefficients; let's call it $C_p(m,n)$ for convenience.
This gives us the following ``sum rule'' for $\delta_p$:
\begin{equation}\tag{$+$}\label{eq+}\delta_p(m+n) = \delta_p(m) + \delta_p(n) + C_p(m,n).
\end{equation}
Products satisfy a rule with a similar flavor:
\begin{equation}\tag{$\times$}\label{eqx}\delta_p(mn) = m^p \delta_p(n) + n^p \delta_p(m) + p \delta_p(m) \delta_p(n).
\end{equation}

The fact that we have rules to break things down into sums and products gives the basis for another comparison with old-fashioned $\frac{d}{dx}$:

\begin{comp}[Sum and product rules]{\ }\label{comp2}
\begin{itemize}
\item For polynomials $f(x), g(x)$, one can compute each of $\frac{d}{dx}(f+g)$ and $\frac{d}{dx}(fg)$ as a fixed polynomial expression in the inputs $f,g,\frac{d}{dx}(f),\frac{d}{dx}(g)$, namely $\frac{d}{dx}(f+g)= \frac{d}{dx}(f)+\frac{d}{dx}(g)$ and $\frac{d}{dx}(fg) = f \frac{d}{dx}(g) + g \frac{d}{dx}(f)$.
\item For integers $m,n$, one can compute each of ${\delta_p(m+n)}$ and $\delta_p(mn)$ as a fixed polynomial expression in the inputs $m,n,\delta_p(m),\delta_p(n)$, namely \eqref{eq+} and \eqref{eqx}.
\end{itemize}
\end{comp}

We might pause to ask whether we could have hoped for a simpler way to differentiate by 13. If we want Comparsion~\ref{comp2} to hold, then the following theorem of Buium provides a definitive answer.

\begin{thm}[{Buium \cite{MR1482984}}] Any function ${\delta:\Z\to \Z}$ that satisfies 
\begin{itemize}
\item a sum rule $\delta(m+n)=S(m,n,\delta(m),\delta(n))$ for some polynomial $S$ with integer coefficients
\item a product rule $\delta(mn)=P(m,n,\delta(m),\delta(n))$ for some polynomial $P$ with integer coefficients
\end{itemize}
is of the form
\[ \delta(n) = \pm \frac{n-n^{p^e}}{p} + f(n)\]
for some prime integer $p$, positive integer $e$, and polynomial $f$ with integer coefficients.
\end{thm}
That is, any function satisfying a sum rule and a product rule is a mild variation on a $p$-derivation.

With the properties of $p$-derivations we have so far, we can recreate analogues of some familiar aspects of calculus.
For example, from the product rule \eqref{eqx} and a straightforward induction, we obtain a power rule:
\[ \delta_p(n^a) = \sum_{i=1}^a \binom{a}{i} p^{i-1} \delta_p(n)^{i} n^{(a-i)p}.\]

Note that the $i=1$ term in the sum above, $a n^{(a-1)p} \delta_p(n)$, looks a bit like the power rule for usual derivatives. If we allow ourselves to extend $\delta_p$ to a map on $\Q$, then we get an analogue of the quotient rule:
\[ \delta_p\left(\frac{m}{n}\right) = \frac{n^p\delta_p(m) -m^p \delta_p(n)}{n^{2p}+pn^p\delta_p(n)}.\]

Of the main cast of characters in a first class on derivatives, perhaps the most conspicuous one missing at this point is the chain rule. Since there is no way to compose a number with a number, we will need a notion of $p$-derivations for functions to state a sensible analogue of the chain rule.

\section*{$p$-derivations for general commutative rings} One can define $p$-derivations for commutative rings with $1$.
\begin{defn} Let $R$ be a commutative ring with $1$ and $p$ a prime integer.
A \emph{$p$-derivation} on $R$ is a function $\delta:R\to R$ such that $\delta(0)=\delta(1)=0$ and $\delta$ satisfies the sum rule $(+)$ and the product rule $(\times)$ above; i.e., for all $r,s\in R$,
\begin{equation}\tag{$+$}\label{eq+2}\delta_p(r+s) = \delta_p(r) + \delta_p(s) - \sum_{i=1}^{p-1} \frac{\binom{p}{i}}{p} r^i s^{p-i}
\end{equation}
and
\begin{equation}\tag{$\times$}\label{eqx2}\delta_p(rs) = r^p \delta_p(s) + r^p \delta_p(s) + p \delta_p(r) \delta_p(s).
\end{equation}
\end{defn}

Evidently, the functions $\delta_p$ we defined on $\Z$ above are $p$-derivations. In fact, for a fixed $p$, a simple induction and the sum rule show that for any $p$-derivation $\delta$ on a ring $R$ and any $n$ in the prime subring (image of $\Z$) of $R$, $\delta(n)=\delta_p(n)$.

The other basic example is as follows. Take the ring of polynomials in $n$ variables with integer coefficients, $R=\Z[x_1,\dots,x_n]$. For any polynomial $f(x_1,\dots,x_n)$, we can consider its $p$th power $ f(x_1,\dots,x_n)^p$, or we can plug in $p$th powers of the variables as inputs to get $f(x_1^p,\dots,x_n^p)$.
These are different, but they agree modulo $p$ as a consequence of the ``Freshman's Dream''. Namely, in the quotient ring $R/pR\cong \Z/p\Z[x_1,\dots,x_n]$,
\[(f+g)^p=f^p + \sum_{i=1}^{p-1} \binom{p}{i} f^i g^{p-i} + g^p = f^p + g^p,\]
since each $ \binom{p}{i}$ is a multiple of $p$, and
\[ (fg)^p = f^p g^p\] as a consequence of commutativity, so the map ${f\mapsto f^p}$ is a ring homomorphism in $R/pR$, called the \emph{Frobenius map}. Thus, in $R/pR$, taking $p$th powers before doing polynomial operations is just as good as after.
So, back in $R$ we \emph{can} divide the difference by $p$, and we \emph{will}! Namely, we can define the function
\[ \delta(f(x_1,\dots,x_n))= \frac{f(x_1^p,\dots,x_n^p) - f(x_1,\dots,x_n)^p}{p},\]
and this function is a $p$-derivation.  Just so we can refer to this function later, let's call this the \emph{standard} $p$-derivation on $\Z[x_1,\dots,x_n]$ and denote it by $\delta_{\st,p}$ (though this notation is not at all standard).

For example,
\[\begin{aligned} \delta_{\st,2}(x^3 + 5x) &= \frac{(x^2)^3 + 5(x^2) - (x^3+5x)^2}{2} \\&= -5x^4 - 10x^2.\end{aligned}\]

As this operator $\delta_{\st,p}$ measures the failure of the Freshman's Dream, one might think of this as a Freshman's Nightmare. In fact, in large generality, $p$-derivations all arise from some freshman's nightmare. Let's make this precise. Given a ring $R$, we say that a map $\Phi:R\to R$ is a \emph{lift of Frobenius} if it is a ring homomorphism and the induced map from $R/pR \to R/pR$ is just the Frobenius map, i.e., $\Phi(r) \equiv r^p \mod pR$ for all $R$. Given a $p$-derivation $\delta:R\to R$, the map  $\Phi:R\to R$ given by \[\Phi(r)=r^p + p \delta(r)\] is a lift of Frobenius. Indeed, the congruence condition is automatic, and the sum rule and product rule on $\delta$ translate exactly to the conditions that $\Phi$ respects addition and multiplication. Conversely, if $p$ is a nonzerodivisor on $R$, and $\Phi$ is a lift of Frobenius, then the map $\delta(r) = \frac{\Phi(r) - r^p}{p}$ is a $p$-derivation: the freshman's nightmare associated to the lift of Frobenius $\Phi$.

It is worth noting that not every ring admits a $p$-derivation. For a quick example, no ring $R$ of characteristic $p$ admits a $p$-derivation, since we would have 
\[0=\delta(0)=\delta(p)=\delta_p(p) = 1-p^{p-1} = 1\] in $R$. Much more subtle obstructions exist, and it is an interesting question to determine which rings admit $p$-derivations; see \cite{MR4269423} for some recent work on related questions. 

Note that the power rule from before follows for any $p$-derivation on any ring, since we just used the product rule to see it. The order-decreasing property holds in general, too, at least if $p$ is a nonzerodivisor on $R$---this follows from writing $s=p^a r$ with $p\nmid r$ and applying the product rule:
\[ \begin{aligned} \delta&(p^a r) = p^{ap} \delta(r) + r^p \delta(p^a) + p \delta(p^a) \delta(r) \\&=  p^{ap} \delta(r) + 
(p^{a-1} - p^{ap-1}) r^p + p (p^{a-1} - p^{ap-1}) \delta(r) \\&=
p^{a-1} r^p + p^a ( \delta(r) - p^{a(p-1)-1} r^p - p^{a(p-1)}\delta(r)).
\end{aligned}\]

Let's wrap up our cliffhanger from the previous section. Now that we have $p$-derivations of polynomials, we have the ingredients needed for a chain rule: given a polynomial $f(x)$ and a number $n$, we will think of the number $f(n)$ as the composition of the function $f$ and the number $n$, and we can try to compare $\delta_p(f(n))$ with $\delta_\st(f)$ and $\delta_p(n)$. Here's the chain rule:
\[ \delta_p(f(n)) = \delta_{\st,p}(f)(n) + \sum_{j=1}^{\deg(f)} {p^{j-1}} \frac{d^j f}{j! \,dx^j}(n^p) \delta_p(n)^j.\]
This is a bit more complicated than the original, but let's notice in passing that the $j=1$ term in the sum, $\frac{df}{dx}(n^p) \delta_p(n)$, looks pretty close to the classic chain rule, besides the $p$th power on $n$. The curious reader is encouraged\footnote{For a hint, consider the lift of Frobenius on $\Z[x]$ that sends $x\mapsto x^p + p \delta_p(n)$, and use Taylor expansion to rewrite the associated $p$-derivation in terms of the standard $p$-derivation and derivatives of $f$.} to prove the formula above.

We have collected a decent set of analogues for the basics of differential calculus for $p$-derivations. One can ask how far this story goes, and the short answer is very far. Buium has developed extensive theories of arithmetic differential equations and arithmetic differential geometry, building analogues of the classical (nonarithmetic) versions of these respective theories with $p$-derivations playing the role of usual derivatives. The reader is encouraged to check out \cite{MR2166202,MR3643159} to learn more about these beautiful theories, though our story now diverges from these. Instead, we will turn our attention towards using $p$-derivations to give some algebraic results with geometric flavors.

\section*{A Jacobian criterion}
\begin{figure}[htb]
\begin{center}
\includegraphics[scale=.4]{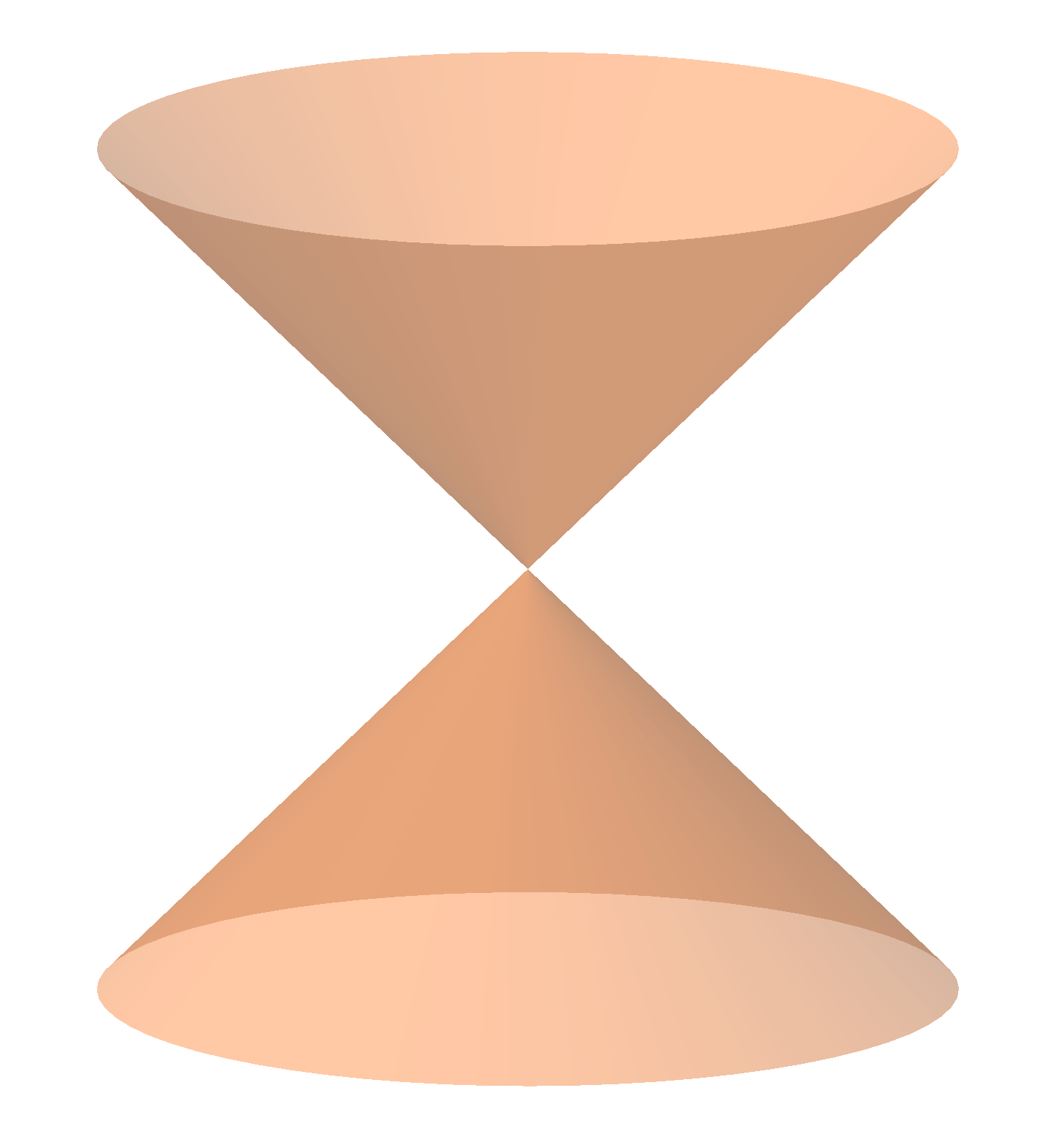}
\end{center}
\caption{ The cone of solutions of $z^2-x^2-y^2=0$. The origin is the unique singular point.}
\label{fig1}
\end{figure}

One natural geometric consideration is whether, and where, a shape has \emph{singularities}: points that locally fail to look flat, due to some sort of crossing or crinkle (or some harder to imagine higher-dimensional analogue or a crossing or crinkle). For example, the double cone cut out by the equation $z^2-x^2-y^2=0$ has a singularity at the origin where the two cones meet, but any other point on the cone is not a singularity, see Figure~\ref{fig1}.

We are going to consider shapes like this that are cut out by polynomial equations, though to state the classical Jacobian criterion, we will consider their solution sets over the complex numbers.

Since it is difficult to envision higher dimensional shapes (and impossible to envision what we're doing next!), it will be useful to give a somewhat more algebraic heuristic definition of singularity. We will say that a point $x$ is a \emph{nonsingular} point in $X$ if within $X$ one can locally cut out $x$ by exactly \mbox{$d=\dim(X)$-many} equations without taking roots, and \emph{singular} otherwise. For example, the point $(1,0,1)$ in the cone is nonsingular, and I claim that the two equations $y=0$, $z-1=0$ ``work'' for our definition: with these two equations and the equation for $X$, we get 
\[ y=z-1=z^2-x^2-y^2=0.\]
Substituting in, we get $0 = x^2-1 = (x-1)(x+1)$, and ``near $(1,0,1)$'', $x+1$ is nonzero, so we can divide out and get $x-1=0$, so $x=1,y=0, z=1$.
On the other hand, $(0,0,0)$ is singular, and the two equations $y=0$, $z=0$ don't ``work'' for our definition: we have 
\[ y=z=z^2-x^2-y^2=0,\]
so $x^2=0$, but we need to take a root to get $x=0$.

The classical Jacobian criterion gives a recipe for the locus of all singularities of a shape cut out by complex polynomials.

\begin{thm}[Jacobian criterion]
Let $X\subseteq \C^n$ be the solution set of the system of polynomial equations
\[ f_1 = \cdots= f_m=0.\]
If the dimension of $X$ is $d=n-h$, and $f_1,\dots,f_m$ generate a prime ideal\footnote{We recall that an ideal  $I$ is prime if it is proper and $gh\in I$ implies $g\in I$ or $h\in I$. Experts will recognize this condition as overkill, but something needs to be done to avoid examples like $(x-y)^2=0$ whose solution set $X$ is the (complex) line $\{(a,a)\}$, which is nonsingular, but $\frac{\partial f}{\partial x} = \frac{\partial f}{\partial y} =0$ for every point in $X$, or $x^2-x=xyz=0$ whose solution set $Y$ is two-dimensional but has a crossing singularity at $(0,0,1)$ that is not a solution of the $1\times 1$ minors of the Jacobian matrix.}  in the polynomial ring $\C[x_1,\dots,x_n]$, then the set of singular points is the solution set within $X$ of the system of polynomial equations
\[\text{all $h \times h$ minors of} \ \begin{bmatrix} \frac{\partial f_1}{\partial x_1} & \cdots & \frac{\partial f_1}{\partial x_n} \\ 
\vdots &  \ddots & \vdots \\
 \frac{\partial f_m}{\partial x_1} & \cdots & \frac{\partial f_m}{\partial x_n}\end{bmatrix} =0.\]
 In particular, if $f$ is irreducible, the set of singular points of the solution set $f=0$ is the solution set of
 \[ f= \frac{\partial f}{\partial x_1}= \cdots = \frac{\partial f}{\partial x_n} =0.\]
 \end{thm}
 
 For example, for the Whitney umbrella cut out by the polynomial $f=x^2-y^2z$, the singular locus is cut out by the system
 \[ x^2-y^2z = 2x = 2yz = y^2 =0,\]
 which simplifies to $x=y=0$; the $z$-axis is where the shape crosses itself.
\begin{figure}[htb]
\begin{center}
\includegraphics[scale=.2]{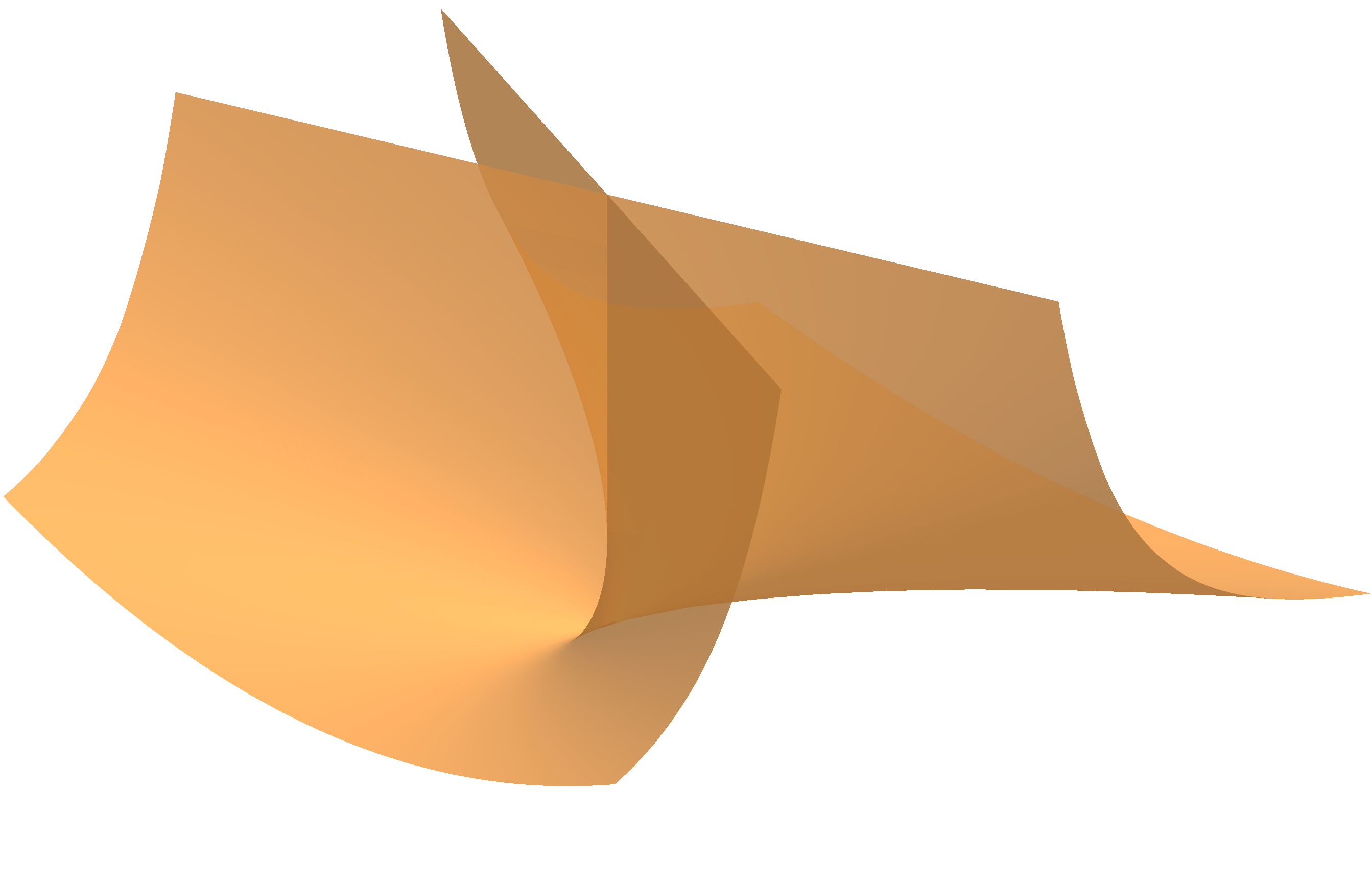}
\end{center}
\caption{ The Whitney umbrella of solutions of ${x^2-y^2z=0}$. The singular locus consists of the line where it crosses itself.}
\end{figure}

The notion of (non)singularity in geometry is generalized in algebra by the notion of regular ring. For starters, prime ideals in algebra play the role of points in geometry: this is motivated by Hilbert's Nullstellensatz, which says that for a quotient ring of the form\footnote{For a set of elements $a_1,\dots,a_t$ in a ring $R$, we use the notation \[(a_1,\dots,a_t)\colonequals \{ r_1 a_1 +  \cdots + r_t a_t \ | \ r_i\in R\}\] for the \emph{ideal generated by $\{a_1,\dots,a_t\}$}.}
\[ \frac{\C[x_1,\dots,x_n]}{(f_1,\dots,f_m)},\]
every maximal ideal
is of the form \[\fm_a=(x_1-a_1,\dots,x_n-a_n)\] for some $a= (a_1,\dots,a_n)$ solution to \[f_1(a) = \cdots=f_m(a)=0;\]
including all prime ideals leads to a better theory for general rings. Then we say that a prime ideal in a ring $R$ is \emph{nonsingular} or \emph{regular} at a prime ideal $\fq$ if $\fq$ can be generated ``locally'' by $h$ equations, where $h$ is the codimension of $\fq$ (how much $\fq$ cuts down the dimension of $R$), and ``locally'' means\footnote{Precisely, ``locally'' means that we work in the localization $R_{\fq}$, and codimension refers to the height of the ideal $\fq$.} that one can divide by elements outside of $\fq$. In the motivating geometric situation where $X$ is the solution set of $f_1=\cdots=f_m=0$ over $\C$, the point $a\in X$ is nonsingular if and only if $R=\C[x_1,\dots,x_n]/(f_1,\dots,f_m)$ is nonsingular at the maximal ideal $\fm_a$. A prime ideal is \emph{singular} if it is not nonsingular.

Intuitively, when working over $\Z$ rather than over $\C$ or a field, in addition to the geometric dimensions, there is an arithmetic dimension that corresponds to the prime integers $p$  in $\Z$. To detect singularity, it suffices to include $p$-derivations as a substitute for derivatives in the $p$-direction! The following is a special case of a result independently obtained by Saito \cite{MR4412577} and Hochster and the author \cite{HJ}.

\begin{thm}[Saito, Hochster-Jeffries]
Let ${R=\displaystyle \frac{\Z[x_1,\dots,x_n]}{(f_1,\dots,f_m)}}$ for some prime ideal $(f_1,\dots, f_m)$ of $\Z[x_1,\dots,x_n]$. Then the set of singular prime ideals $\fq$ of $R$ containing a fixed prime integer $p$ is exactly
the set of prime ideals containing $p$ and
\[\substack{ \text{all $h \times h$}\\ \text{ minors of}} \ \begin{bmatrix} \delta_{\st,p}(f_1) & (\frac{\partial f_1}{\partial x_1})^p & \cdots & (\frac{\partial f_1}{\partial x_n})^p \\ 
\vdots  & &  \ddots & \vdots \\
 \delta_{\st,p}(f_m) & (\frac{\partial f_m}{\partial x_1})^p & \cdots & (\frac{\partial f_m}{\partial x_n})^p \end{bmatrix},\]
 where $h=n+1-\dim(R)$.
 In particular, if $f$ is irreducible, the set of singular prime ideals of $R=\frac{\Z[x_1,\dots,x_n]}{(f)}$ containing $p$ is exactly the set of primes containing $p$ and
 \[ \delta_{\st,p}(f)\ ,\ \frac{\partial f}{\partial x_1}\ ,\ \ldots \ ,\ \frac{\partial f}{\partial x_n} .\]
\end{thm}

\begin{ex} Let $n$ be a squarefree integer (excluding $0$ and $1$), $q$ a prime number, and consider the ring $R=\Z[\sqrt[q]{n}]$. We claim that this admits a singular prime ideal if and only if $\delta_q(-n)$ is a multiple of $q$. 
Think\footnote{This is where excluding $0$ and $1$ is necessary.} of $R$ as $\Z[x]/(x^q-n)$. For a prime number $p$, the singular prime ideals $\fa$ containing $p$ are those that contain
\[ p , \delta_{\st,p}(x^q-n), \frac{d}{dx}(x^q-n).\]
Using the sum rule for $\delta_{\st,p}$ and the defining equation for $R$, we have 
\[ \delta_{\st,p}(x^q-n) = \begin{cases} \delta_p(-n)  &\text{ if $p\neq 2$}\\  \delta_p(-n) - nx^q  &\text{ if $p= 2$}\end{cases}\]  in $R$. 

For $p\neq q$, from $\frac{d}{dx}(x^q-n) = qx^{q-1}$, we get that $\fa$ must contain $x$, and from the defining equation, $n$ as well. But if the integers $p,n,\delta_p(n)$ are in a proper ideal $\fa$, since $p$ is a prime number, we must have $p\mid n$ and $p\mid \delta_p(n)$, since $1$ would be a linear combination of these numbers otherwise. 
By the Order Decreasing Property, $p^2 \mid n$, contradicting that $n$ is squarefree. So, there are no singular prime ideals containing ${p\neq q}$.

For $p=q\neq 2$ since $\frac{d}{dx}(x^q-n)$ is a multiple of $q$, using the simplification above, a prime ideal containing $q$ is singular if and only if it contains 
\[ q , \delta_q(-n).\]
But if $q$ and $\delta_q(-n)$ are in $\fa$, then $q \mid \delta_q(-n)$; a singular prime ideal then occurs if and only if this happens. The analysis for $p=q=2$ is similar (cf. \cite{HJ}).

This has a consequence for a familiar object in elementary number theory. By some standard results\footnote{Namely, $\Z[\sqrt[q]{n}]$ is a ring of integers if and only if it is integrally closed in its fraction field. Since this ring is an integral extension of $\Z$ generated by one element as an algebra, it is a one dimensional Noetherian domain. Such a ring is integrally closed in its fraction field if and only if it is nonsingular.} in commutative algebra, the ring $\Z[\sqrt[q]{n}]$ is a ring of algebraic integers (for its fraction field $\Q(\sqrt[q]{n})$) if and only if it has no singular prime ideals. Thus, we conclude that $\Z[\sqrt[q]{n}]$ is a ring of integers if and only if $q\nmid \delta_q(-n)$. In particular, for $q=2$, using our earlier observation on when $\delta_2(n)$ is odd or even, we recover the fact that for $n$ squarefree, $\Z[\sqrt{n}]$ is a ring of integers if and only if $n\not\equiv 3 \mod 4$.
\end{ex}

\section*{A Zariski-Nagata theorem for symbolic powers}

Let's recall another classical theorem relating algebra and geometry. To state it, we need the notion of symbolic power of a prime ideal. Over a polynomial ring over a field, or more generally, over a commutative Noetherian ring, any ideal can be written as an intersection 
\begin{equation}\label{eq:pd}\tag{3}
I = \fq_1 \cap \cdots \cap \fq_t
\end{equation}
where the $\fq_i$ are \emph{primary} ideals: ideals with the property that $xy\in \fq \Rightarrow x \in \fq \text{ or } y^n\in \fq$ for some $n$. Such an expression is called a \emph{primary decomposition}. The existence of primary decomposition is a famous result of Lasker in the polynomial case and Noether in the Noetherian case. This can be thought of as a generalization of the Fundamental Theorem of Arithmetic: in $\Z$, the primary ideals are just the ideals generated by powers of primes, and a prime factorization
\[ n= p_1^{e_1} \cdots p_n^{e_n} \]
corresponds to writing $(n)$ as an intersection of primary ideals:
\[ (n) = (p_1^{e_1}) \cap  \cdots \cap (p_n^{e_n}).\]
There are two important differences with the Fundamental Theorem of Arithmetic, though:
\begin{itemize}
\item primary ideals are not powers of prime ideals in general, nor are powers of prime ideals always primary, and
\item the collection of primary ideals $\fq_i$ appearing in the decomposition \eqref{eq:pd} is not unique, but if the decomposition satisfies a simple irredundancy hypothesis, then the components whose radical does not contain any other component's radical are uniquely determined.
\end{itemize}

In particular, if $\fq$ is a prime ideal and $n>1$, then according to the first point above, $\fq^n$ may admit an interesting primary decomposition, and as a consequence of the second point, the component with radical $\fq$ is the same in any (irredundant) primary decomposition. This is called the $n$th \emph{symbolic power} of $\fq$, denoted $\fq^{(n)}$.

Symbolic powers arise in various contexts in algebra and geometry. For example, they arise naturally in interpolation questions, they play a key role in the proofs of various classical theorems such as Krull's Principal Ideal Theorem, and they have enjoyed a resurgence of interest in combinatorics in connection with the Packing Problem of Conforti and Cornu\'ejols. The interest reader is recommended to read the survey \cite{MR3779569}.

A classical pair of theorems of Zariski and Nagata gives a geometric description of the symbolic power of an ideal in a polynomial ring over $\C$. The result has various statements; we will give a differential statement.

\begin{thm}[Zariski-Nagata Theorem] Let ${X\subseteq \C^n}$ be the solution set of the system of polynomial equations
\[ f_1 = \cdots= f_m=0.\]
Suppose that $f_1,\dots,f_m$ generate a prime ideal~$\fq$. Then $\fq^{(r)}$  is exactly the set of polynomials ${f\in \C[x_1,\dots,x_n]}$ such that
\[ \frac{\partial^{a_1+\cdots+a_n} f}{\partial x_1^{a_1} \cdots \partial x_n^{a_n}}  \Big|_{X} \equiv 0 \text{ for all } a_1+\cdots+a_n < r.\]
\end{thm}

The same characterization is doomed to fail in $\Z[x]$: for example, $\fq=(2)$ is a prime ideal with $\fq^{(2)} = (4)$; in particular, $2\notin (4)$. But $2$ satisfies the derivative condition corresponding to the right hand side above: taking $a_1=0$, we have $2\in (2)$, and taking $a_1=1$, we have $\frac{\partial 2}{\partial x} = 0 \in (2)$. 
 
If you've been paying attention so far, you should be able to name the missing ingredient.  Indeed, $\delta_{\st,2}(2)=-1 \notin (2)$, so allowing partial derivatives and a $2$-derivation is enough to take this element $2$ that isn't in $\fq^{(2)}$ out of $\fq$, suggesting a way to characterize $f\in \fq^{(2)}$ in terms of derivatives (including our ``derivative by 2'').

 In fact, this works in general. The following analogue of the Zariski-Nagata Theorem is a special case of a result of De Stefani, Grifo, and the author \cite{MR4080246}.

\begin{thm}[De Stefani-Grifo-Jeffries]
Consider the ring $\Z[x_1,\dots,x_n]$ and let $\fq=(f_1,\dots,f_m)$ be a prime ideal. Suppose that $\fq$ contains the prime integer $p$. Then $\fq^{(r)}$ is exactly the set of polynomials $f\in \Z[x_1,\dots,x_n]$ such that
\[ \delta_{\st,p}^{a_0}\left(\frac{\partial^{a_1+\cdots+a_n} f}{\partial x_1^{a_1} \cdots \partial x_n^{a_n}} \right)  \in \fq \text{ for all } a_0+a_1+\cdots+a_n < r.\]
\end{thm}

\section*{Other applications}
 $p$-derivations have appeared in a range of sophisticated applications to number theory and arithmetic geometry. We briefly list a few of these, and encourage the reader to explore further.

 \subsection*{Effective bounds on rational points and $p$-jet spaces} The motivation for Buium's original work on $p$-derivations was to give bound the number of points on rational points on curves. For a complex curve $X$ defined over $\overline{\Q}$, there is a natural map from $X$ to an algebraic group $A$, the Jacobian of $X$; the main result of Buium \cite{MR1387233} gives an effective bound, depending only on the genus of $X$ and the smallest prime of a good reduction for $X$, on the number of points in $X$ that map to torsion elements of $A$.
 
To establish these results, Buium constructs \emph{$p$-jet spaces}. In differential geometry, the \emph{jet} of a function $f$ of order $k$ at a point is the data of all of the values of the derivatives of $f$ up to order $k$ at that point; the \emph{jet space} of order $k$ of a manifold $X$ is a manifold whose points correspond to jets on $X$. Buium's $p$-jet spaces are analogues of jet spaces obtained by replacing usual derivatives with $p$-derivations. The result mentioned above is then obtained by intersection theory on $p$-jet spaces of curves.

As mentioned earlier, Buium has developed an extensive theory of arithmetic differential geometry in analogy with classical differential geometry, for which $p$-jet spaces form the starting point. We refer the reader to \cite{MR2166202} to learn more.
  
 \subsection*{Relationship with Witt vectors} The Witt vectors are a construction in number theory that generalizes the relationship between the prime field $\F_p$ and the corresponding ring of $p$-adic integers $\Z_p$. To be precise, there is a functor $W$ from the category of rings to itself, called the functor of ($p$-typical) Witt vectors, that maps $\F_p$ to $\Z_p$. More generally, for any perfect field $F$ of characteristic $p$, $W(F)$ is a local ring with maximal ideal $(p)$, and $W(F)/pW(F) \cong F$; in this way, one can think of $W$ as a generalization of a construction of $\Z_p$ from $\F_p$.
 
$p$-derivations have many interesting connections with Witt vectors; indeed, they first arise in Joyal's work to study Witt vectors. Namely, Joyal shows that the forgetful functor from the category of rings with a $p$-derivation to the category of rings is left adjoint to the Witt vector functor.

 \subsection*{Philosophy of the field with one element} Various formulas for enumerating basic objects over a finite field $\F_q$ limit to combinatorially meaningful quantities as $q\to 1$; for example, 
 \[ \lim_{q\to 1} \frac{| \mathrm{GL}_n(\F_q) |}{|(\F_q ^\times)^n|} = n!,\]
 where $\mathrm{GL}_n(\F_q)$ is the collection of linear automorphisms of $\F_q^n$, the vector space of $n$-tuples over $\F_q$, and $n!$ counts the number of permutations of $n$ elements.
 Partially motivated by this phenomenon, and partially motivated by transferring results over finite fields to other settings, there is a program of inventing a notion of algebraic geometry over a ``field with one element'': the mythical field with one element is not literally a field with one element, which would contradict the definition of field, but something with a different structure than a field that admits some sort of geometry analogous to what one would expect over a field with one element for quantitative or various other reasons. On the algebraic side, the field with one element can be thought of as a ``deeper base ring'' than $\Z$ (though not necessarily itself a ring!).
 
There are many approaches in the literature to implementing the philosophy of the field with one element. Most relevant for this article is the theory established by Borger \cite{Bor} as well as the closely related approach of Buium. Roughly speaking, Borger proposes that a model for the field with one element should be $\Z$ equipped with the collection of all its $p$-derivations $\delta_p$. In particular, in Borger's model, $\mathrm{GL}_n$ ``over the field with one element'' is the symmetric group on $n$ letters, aligning with the numerical coincidence noted in the previous paragraph.
 
 \subsection*{Unifying cohomology theories} Recent work of Bhatt and Scholze \cite{MR4502597} has employed $p$-derivations to relate various $p$-adic cohomology theories (\'etale, de Rham, crystalline). The key ingredient is the notion of a \emph{prism}, which is a ring equipped with a $p$-derivation subject to some conditions. We refer the reader to  \cite{MR4502597}  to learn more.
 
\section*{Acknowledgments}

The author is grateful to Alessandro De Stefani, Elo\'isa Grifo, Claudia Miller, Steven Sam, and the referee for many helpful comments on this article.

\bibliographystyle{foo}
\bibliography{ExampleRefs}

\end{document}